# Automatic Analysis, Decomposition and Parallel Optimization of Large Homogeneous Networks

*[1] D.Yu. Ignatov <ignatov.dmitry@huawei.com>*
*[1] A.N. Filippov <filippov.alexander@huawei.com>*
*[2] A.D. Ignatov <ihnatova@student.ethz.ch>*
*[1] X. Zhang <zhangxuecang@huawei.com>*
*[1] Russian Research Center, Huawei Technologies, 1, Altufevskoe sh., Moscow, 127106, Russia*
*[2] Swiss Federal Institute of Technology in Zurich, Department of Computer Science, 101, Rämistrasse, Zurich, 8092, Switzerland*

**Abstract**. The life of the modern world essentially depends on the work of the large artificial homogeneous networks, such as wired and wireless communication systems, networks of roads and pipelines. The support of their effective continuous functioning requires automatic screening and permanent optimization with processing of the huge amount of data by high-performance distributed systems. We propose new meta-algorithm of large homogeneous network analysis, its decomposition into alternative sets of loosely connected subnets, and parallel optimization of the most independent elements. This algorithm is based on a network-specific correlation function, Simulated Annealing technique, and is adapted to work in the computer cluster. On the example of large wireless network, we show that proposed algorithm essentially increases speed of parallel optimization. The elaborated general approach can be used for analysis and optimization of the wide range of networks, including such specific types as artificial neural networks or organized in networks physiological systems of living organisms.

## 1. Introduction

There are many large homogeneous networks consisting of the elements of the same type, which influence our everyday life, such as wired or wireless networks, network of switches in datacenter or artificial neural networks. In addition to mentioned information systems, there are networks transferring physical objects, such as a road network, different types of pipe networks or even organized in networks physiological systems of living organisms. The support of the effective continuous functioning of such networks requires their permanent screening and optimization. It is important that optimization of artificial networks begins on the stage of their planning, continues during all time of their usage and includes balancing of activity of network elements. Thus, whenever we begin an optimization, the network is already partially optimized and its elements have relatively similar level of activity.

An optimization of networks can be described as a maximization of objective function, which simulates activity of network. Objective function takes as input network parameters and provides estimation of performance and/or quality of service provided by network. The homogeneous networks have useful properties – since they consist of uniform elements, the network or its part can be optimized with the usage of the same objective function. Therefore, network can be decomposed into the relatively independent subnets, which are optimized in parallel processes. The quality objective function can be specified, for example, as an average level of radio signal in the area covered by the wireless network, or as an average speed of traffic in switches of datacenter or on crossroads of road network. In general, for homogeneous network decomposed into non-overlapping subnets, the quality objective function can be represented as an average value of objective functions of subnets:

$$Q\left(\vec{w}\right) = \frac{1}{n}\sum_{i=1}^{n} Q\left(\vec{w_i}\right),$$

where Q – quality of service provided by network, $\vec{w}$ – vector of network parameters, n – number of non-overlapping subnets, $\vec{w_i}$ – vector of parameters of $i$-th subnet.

The decomposition of complex systems before their optimization is an old idea raised with development of large-scale algorithms. Beginning from early works [1, 2] until now [3-5] the main decomposition approach is the detection of the mostly independent subnets with minimal strength of interactions between their elements. The network analysis leads to the calculation of strength of these interactions, and then decomposition is performed according to specific for network criterion. Typically, decomposition of network falls under the category of NP-hard problems, which are solved with graph partitioning algorithms. Obtained subnets can be optimized on multiple-core computer or cluster with such general technique as Simulated Annealing [4] or network-specific optimization algorithm [5].

The decomposition allows to do the optimization in parallel as well as to reduce optimization complexity by discarding of weak interactions [4], and, in such way, essentially increases speed of optimization. Remarkable, that decomposition of large network into the weakly interacting parts can reduce the complexity of error surface relief, and in such way decreases probability of stuck in local optimum [6]. Therefore, decomposition of network reduces optimization complexity as well as probability of optimization stuck in local optimum.

The quantity of networks, consisting of uniform elements, rapidly grows in modern industry, especially in the sphere of information technologies. Such growth stimulates the development of new efficient algorithms of network decomposition, which can provide maximal decrease of the optimization complexity and efficient parallelization of optimization process. In this paper, we propose optimization methods, which combine idea of independent optimization with the alternative decomposition approach. They are intended to solve problems with following characteristics: first, optimized network consists of the uniform elements; second, optimization of network is carried out in multi-core and distributed systems.

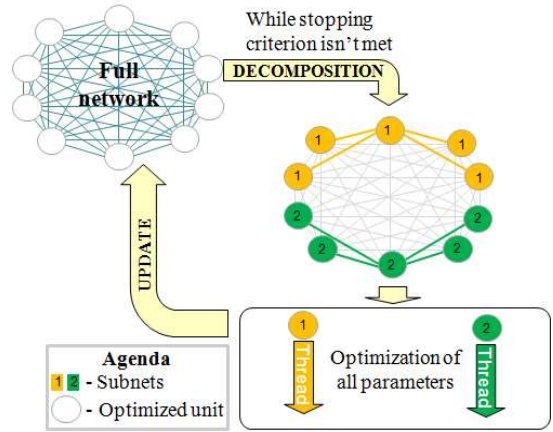

*Fig. 1. Automatic sector planning for parallel optimization of network.*

## 2. Network Decomposition approaches

With the purpose to find solution of the Maximal Independent Set and related problems the big variety of decomposition algorithms are elaborated [1-5]. They break a complex task of network optimization down into subtasks by decomposition of network into one set of relatively independent subnets. In order to perform decomposition, network is firstly represented with the weighted complete graph, where each vertex corresponds to a network element and each edge has weight equals to the rank of correlation between pair of correspondent elements. Hereby the rank of correlation is calculated by network-specific correlation function and represents the strength of elements interactions within network. The weights on edges are used in the process of network decomposition into subnets weakly connected with each other, while data in vertex – for optimization of obtained subnets.

For example, described in invention [3] method of automatic sector planning of network for parallel optimization (Fig. 1) decomposes network into a predefined number of subnets by the rule of minimal sum of the crossing edges weights. After decomposition the subnets are optimized independently with the same objective function. Optimization of subnets is performed for all regulated parameters in distributed-centralized mode.

Main drawback of such approaches is an ignorance of crossing edges, while some of them appear to be significant. Thus, this ignorance negatively affects accuracy and time of optimization process.

## 3. Alternative decomposition with independent optimization

The alternative decomposition (Fig. 2) as well as existing approach (Fig. 1) is based on the representation of network data and structure with the weighted complete graph.

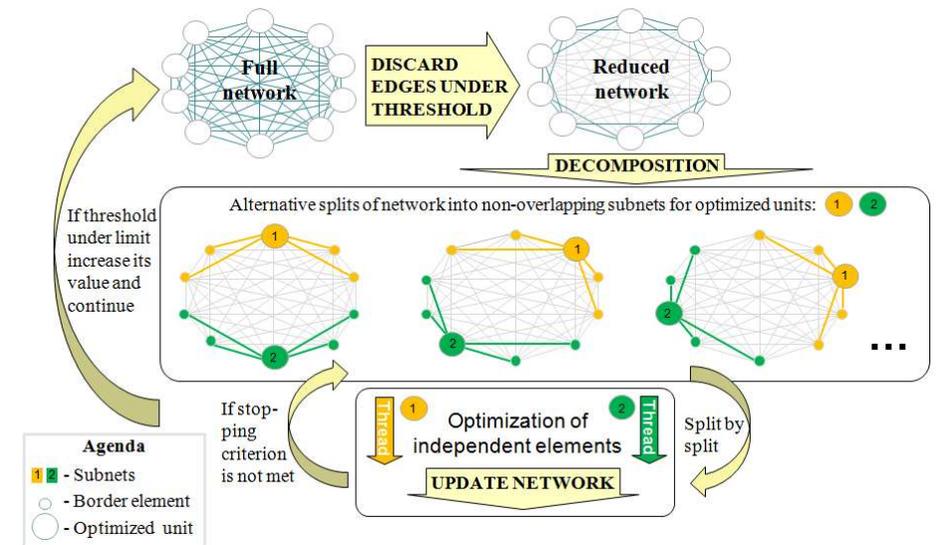

*Fig. 2. Alternative decomposition of network with independent optimization.*

The decomposition consists of two steps:
1. The weak correlations between elements are discarded by removing of edges with weights under network-specific filtering threshold.
2. Reduced network is decomposed into the set of alternative splits containing non-overlapping subnets.

Last step includes two stages (Fig. 3):
  a) Selection of subnet for every optimized unit. This subnet consists of the unit and connected to it elements.
  b) Finding (i.e. brute force search) of combinations of non-overlapping subnets – alternative splits of network. Every split is generated from at least one untapped subnet and covers as much vertices as possible.

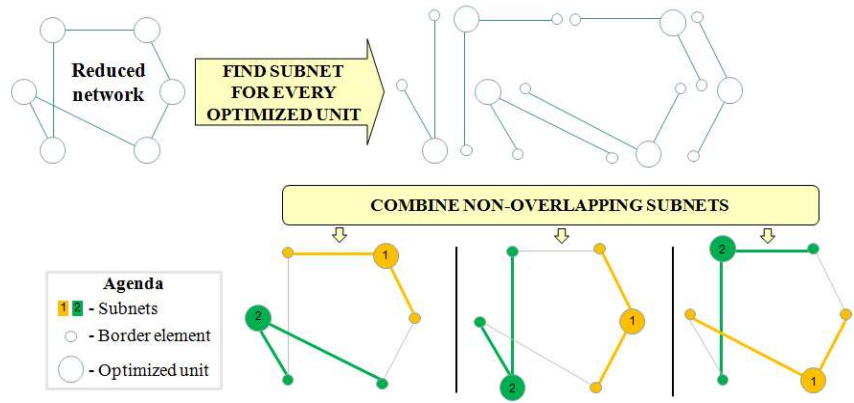

*Fig. 3. Details of alternative decomposition.*

Therefore, in order to reduce complexity of the optimization process the most irrelevant interactions between elements of network are removed. Based on remaining connections, the network is further split into relatively independent subnets. The result of splitting process is a set of alternative splits, where every split consists of all possible non-overlapping subnets within network. The obtained splits are optimized one by one in loop (Fig. 2, bottom). Within each iteration, the following steps are executed:
  1. The next alternative split is selected.
  2. The non-overlapping subnets are apportioned to the cores available on computer/cluster.
  3. The optimization of subnets is performed in parallel processes by optimizing procedure.
  4. The full network is updated with the values of optimized parameters.

Herewith the optimizing procedure is implementation of Simulated Annealing probabilistic technique of global optimization, which randomly searches for optimal solution in the space of all available alternatives [7]. We modified this method by adding new variable – step regulated by precision parameter (P), which restricts the area of neighbors search.

```
    procedure optimize(S0, P) {
      Snew := S0
      step := maxStep * (1 – P)
      Sgen := random neighbor of S0 within step
      T := temperature(1 – P)
      if A(E(S0), E(Sgen), T) ≥ random(0, 1) then Snew := Sgen
      Output: state Snew
    }
```

where $S_0$, $S_{gen}$ and $S_{new}$ – current, generated and new states of subnet, correspondently; maxStep – maximal value of step; temperature – monotonically increasing function mentioned in Simulated Annealing method; E and A – energy and acceptance functions of Simulated Annealing method, correspondently.

This optimizing procedure is performed on alternative splits in cycle as long as better states are founded. After optimizing procedure (Fig. 2, bottom) is finished, we increase filtering threshold value and repeat the cycle: network decomposition and further iterative optimizing procedure. Finally, when the filtering threshold reaches the maximal value, then the optimization process is stopped and optimized parameters are ready to be used for adjustment of physical network.

The effectiveness of this algorithm (Fig. 2) is provided by the parallel optimization of the most independent elements only, whereas in original approach (Fig. 1) all elements of network are optimized at once. The following regulation of optimization precision gives additional increase of optimization speed and quality. We use precision parameter P to control the precision of optimization by regulation of the step of Simulated Annealing algorithm as well as for calculation of threshold for filtering out weak connections between elements. At the beginning of the process, value of P is set to its minimum (i.e. zero) and is progressively increased with predefined constant up to the maximum (i.e. one). Before every decomposition procedure, the value of P is increased

$$P = P + const$$

and filtering threshold is calculated

$$Th = Th_{min} + (Th_{max} - Th_{min}) \times P,$$

where Th – filtering threshold parameter of decomposition process, $Th_{min}$ – network-specific minimal value of Th, $Th_{max}$ – network-specific maximal value of Th, P – precision parameter of optimization process.

Fig. 2 and 3 represent an example of simple network decomposition into set of alternative splits, where optimized unit consists of one element. In the case of large networks, decomposition could be based on the complex optimized units – group of strongly correlated elements, which are optimized together within subnet. The size of optimized unit is determined by the size of network and number of cores available on computer or cluster.

The peculiarities of optimization of alternative splits in distributed system are described bellow (Fig. 4). At the beginning of optimization process the copies of network-specific optimizing procedure are distributed on all nodes. Then the network is decomposed into alternative splits of subnets. Herewith for every split while the number of subnets is lower than the total quantity of cores in computer cluster, subnets are randomly duplicated. When the number of subnets becomes equal to the quantity of cores in cluster, then optimization begins and subnets are sent to every node proportionally to quantity of cores on it. Optimization is carried out independently on every node without data shift. After optimization is finished every node sends back only optimized (i.e. changed) parameters of subnets, and if the alternative solutions exist for some subnet – the best of them is selected. Then network is updated with optimized parameters and the next split is processed.

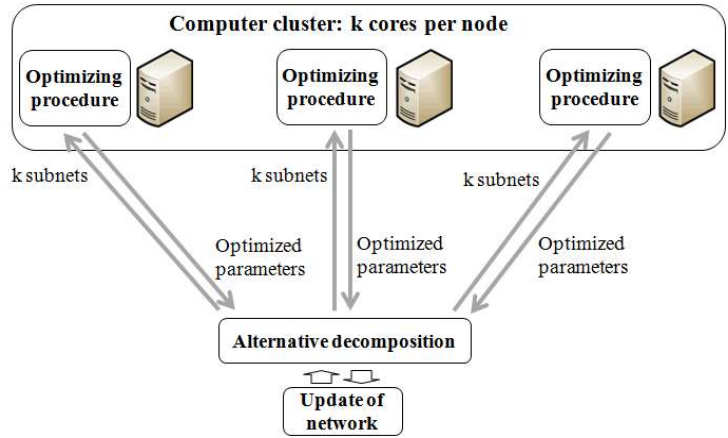

*Fig. 4. Optimization of alternatives on computer cluster.*

Optimization continues as long as objective function produces better result. If optimization of predefined number of splits does not give improvement, the optimization loop is finished, filtering threshold is increased and optimization continues on the next level of splitting. It is important, that in every split the optimized parameters are selected in such way, that they influence activity of full network. Thus, if optimization of predefined number of splits cannot improve network, the probability that the next split will give significant improvement is quite small. The adjustment of physical network is provided continuously – at the end of every optimization loop – or once, when the P value reaches its maximum.

## *4. Evaluation*

For demonstration of the efficiency of proposed approach, the optimization of wireless network is implemented with visualization of quality of radio signal (Fig. 5). The quality is estimated with Signal to Interference plus Noise Ratio (SINR) [8]. An optimizing procedure is represented by modified Simulated Annealing algorithm. Objective function calculates the average value of SINR in covered by network area on the basis of Okumura-Hata, COST-231 and Stanford University Interim radio propagation models [9], which take into account the type of area and radio-frequency diapason. The result of radio propagation model is adjusted according to tilt and azimuth propagation functions [10]. The experimental model of network consists of 100 sites with 3 sector antennas in each of them (total 300 antennas). Every antenna has 4 regulated parameters: power, height, tilt and azimuth (total 1200 regulated parameters). The rank of correlation between pair of antennas is calculated as 1 / (distance between antennas).

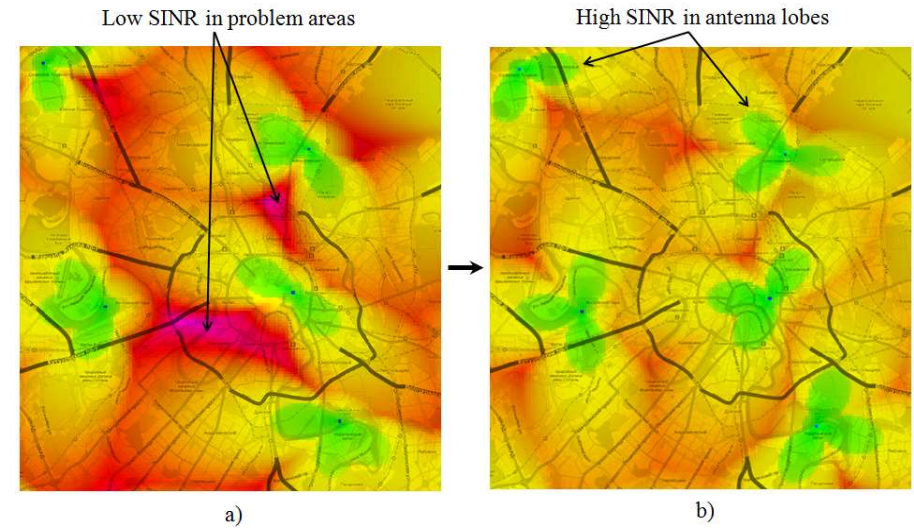

*Fig. 5. Visual representation of Signal to Interference plus Noise Ratio (SINR) in initial (a) and optimized (b) subnet. Crimson and red territories – problem areas with low level of SINR, orange and yellow – sufficient level of SINR, and green – high level of SINR.*

The SINR is automatically registered during optimization and its average values for 30 experiments are represented in logarithmic scale (Fig. 6). We can see, that after 3600 s the values (± 95 % conf. interval) of SINR reaches for sector planning 35.0 ± 1.7 dB and for independent optimization of alternatives 40.1 ± 1.9 dB, and don't change within experimental error after this time. Difference between represented values are significant ($p < 0.01$) and for independent optimization of alternatives 15 % higher, than in case of sector planning. According to approximation curve, the independent optimization of alternatives reaches the value 35.0 dB at ~ 400 s of experiment (Fig. 6, dashed line), so it is ~ 9 times faster, than sector planning algorithm with full optimization. Thus, by the example of large wireless network optimization we show that compare to recently published technique [3] the proposed algorithm (Fig. 2) gives 9 times speed-up or after long time optimization – demonstrates better accuracy.

The main principles, which according to our opinion underlay the speed-up of optimization, are represented on diagram (Fig. 7). According to these principles, optimization process begins with rough search of global optimum in small number of complex subnets with the aim to avoid stuck in local optimum. During optimization the complexity of subnets is progressively decreased without loosing of optimization precision. Such decrease is reasonable, because during optimization process, the amplitude of oscillation of regulated parameters is decreased, and as a result, the strength of close interaction of element within network can essentially prevalent over the strength of distant interactions. According to our experiments (Fig. 6), this strategy gives accurate result of the optimization process.

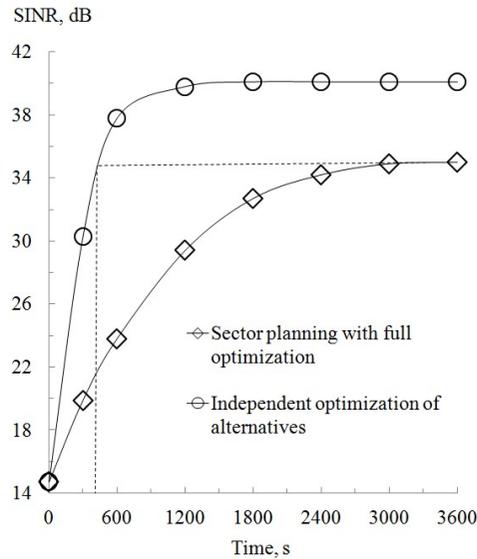

*Fig. 6. Dependence of average value of Signal to Interference plus Noise Ratio (SINR) on time of optimization of wireless network.*

From other side, automatic regulation of quantity of alternative calculations allows us to use all available computational resources for optimization. It is important that at the beginning of optimization process the rough search of global optimum produces essential oscillation of optimized parameter. Therefore, in this situation the alternative calculations are reasonable for selection of the best solution. Oppositely, at the end of the process, the optimization in small subnets produces almost the same results and the quantity of alternative calculations, which are not justifiable any more, are decreased (Fig. 7) to minimum. Thus, proposed algorithm accelerates the optimization process through regulation of precision parameter together with rational usage of all available calculation resources at different time of optimization process.

## 5. Conclusion

This paper presents a new meta-algorithm of network alternative decomposition into the sets of non-overlapping subnets and parallel optimization of the most independent parameters of network on computer or cluster. The optimization is carried out with progressive decreasing of subnet complexity and increasing of optimization precision. Usage of proposed algorithm leads to the following benefits:

- the faster optimization due to reduction of optimization problem complexity and efficient usage of all calculation resources;
- the better precision of optimization due to progressive change of optimization strategy from rough search of optimum at the beginning of optimization process to precise search of optimum at the end of optimization process.

We believe that proposed approaches of independent alternative optimization and dynamic regulation of precision are providing solid basis for the implementation of highly scalable distributed solution for the wide variety of large homogeneous networks.

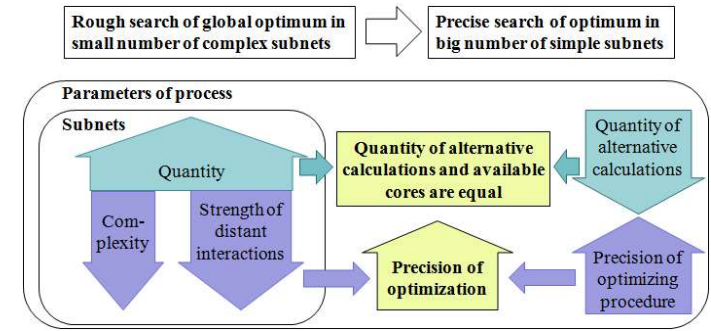

*Fig. 7. Directions of progressive change of parameters of decomposition and optimization with the time. Vertical arrows: up and down – increase and decrease of parameter value, correspondently; horizontal arrows – influence of parameters on integral indices.*

# Автоматический анализ, декомпозиция и параллельная оптимизация больших однородных сетей

*[1] Д.Ю. Игнатов <ignatov.dmitry@huawei.com>*
*[1] А.Н. Филиппов <filippov.alexander@huawei.com>*
*[2] А.Д. Игнатов <ihnatova@student.ethz.ch>*
*[1] С. Чжан <zhangxuecang@huawei.com>*

*[1] Российский научно-исследовании центр, Техкомпания Хуавэй, 127106, Россия, г. Москва, Алтуфьевское шоссе, д. 1.*
*[2] Швейцарская высшая техническая школа Цюриха, отдел компьютерных наук, 8092, Швейцария, г. Цюрих, Рамиштрассе, д. 101.*

**Аннотация**. Жизнь современного мира во многом зависит от функционирования больших однородных сетей, таких как проводные и беспроводные коммуникационные системы, сети дорог и трубопроводов. Поддержание их эффективной работы требует автоматического контроля, постоянной оптимизации, включающей обработку больших объемов данных с использованием высокопроизводительных распределенных систем. Предложен новый мета-алгоритм для анализа больших однородных сетей, их альтернативного разбиения на слабосвязанные подсети и параллельной оптимизации наиболее независимых элементов подсетей. Данный подход основан на специфической для сети корреляционной функции, алгоритме имитации отжига и адаптирован для работы в вычислительном кластере. На примере беспроводной коммуникационной сети показано, что предложенный алгоритм существенно увеличивает скорость многопоточной оптимизации. Разработанный общий подход может быть использован для анализа и оптимизации широкого спектра сетей, включая такие специфические типы как искусственные нейронные сети или организованные в виде сетей физиологические системы живых организмов.

**Ключевые слова:** однородные сети; декомпозиция; оптимизация; распределенные вычисления

## Список литературы